\documentclass[12pt]{article}
\usepackage[utf8]{inputenc}
\usepackage{amssymb,amsmath,amsfonts,eucal,mathrsfs,amsthm} 
\usepackage[colorinlistoftodos]{todonotes}
\usepackage[allcolors=blue]{hyperref}
\usepackage{empheq}
\newtheorem{theorem}{Theorem}

\newtheorem{lemma}[theorem]{Lemma}

\newtheorem{corollary}[theorem]{Corollary}

\newtheorem{remark}[theorem]{Remark}

\theoremstyle{definition}
\newcommand{\R}{\mathbb{R}}
\newcommand{\Q}{\mathbb{Q}}
\newcommand{\Sf}{\mathbb{S}}

\newcommand{\Hy}{\mathbb{H}}

\newcommand{\trace}{\mbox{tr\,}}

\def\<{{\langle}}
\def\>{{\rangle}}

\def\n{\nabla}
\def\d{\partial}
\def\a{\alpha}

\def\vol{{\rm{Vol}}}
\def\be{\begin{equation} }
\def\ee{\end{equation} }

\def\proof{\noindent{\it Proof:  }}
\def\qed{\ifhmode\unskip\nobreak\fi\ifmmode\ifinner
\else\hskip5 pt \fi\fi\hbox{\hskip5 pt \vrule width4 pt
height6 pt  depth1.5 pt \hskip 1pt }}
\makeatletter
\newcommand{\subjclass}[2][]{\let\@oldtitle\@title
\gdef\@title{\@oldtitle\footnotetext{#1 
\emph{Mathematics Subject Classification:} #2}}}
\newcommand{\keywords}[1]{\let\@@oldtitle\@title
\gdef\@title{\@@oldtitle\footnotetext
{\emph{Key words and phrases.} #1.}}}
\makeatother
\begin{document}

\title{Isometric immersions with flat normal bundle
between space forms}
\author{M. Dajczer, C.-R. Onti and Th. Vlachos}
\date{}
\maketitle

\begin{abstract}
We investigate the behavior of the second fundamental form 
of an isometric immersion of a space form with negative curvature 
into a space form so that the extrinsic curvature is negative. 
If the immersion has flat normal bundle, we prove that its second 
fundamental form grows exponentially. 
\end{abstract}

It is a long-standing problem if the complete hyperbolic space 
$\Hy^n$ can be isometrically immersed in the Euclidean space 
$\R^{2n-1}$. In fact, the non-existence of such an immersion 
has been frequently conjectured; see Yau \cite{yau}, Moore 
\cite{moorepro} and Gromov \cite{gromov}. A positive answer
to the conjecture would be a natural generalization to higher 
dimensions of the classical result from $1901$ by Hilbert for the 
hyperbolic plane. On one hand, Cartan \cite{ca} in 1920
showed that $\Hy^n$, $n\geq 3$, cannot be isometrically immersed 
in $\R^{2n-2}$ even locally. On the other hand, he proved that 
there exists an abundance of local isometric immersions of $\Hy^n$ 
into $\R^{2n-1}$ and that these  have all flat normal bundle.

Nikolayevsky \cite{niko} proved that complete non-simply connected
Riemannian manifolds of constant negative  sectional curvature 
cannot be isometrically immersed into Euclidean space with flat 
normal bundle. Let $\Q^m_c$ denote a complete simply connected 
$m$-dimensional Riemannian manifold of constant sectional curvature 
$c$, that is, the Euclidean space $\R^m$, the Euclidean sphere 
$\Sf_c^m$ or the hyperbolic space $\Hy_c^m$ according to whether 
$c=0, c>0$ or $c<0$, respectively. It was observed in \cite{dt} 
that the proof by Nikolayevsky gives, in fact, the following
slightly more general result:
\vspace{1ex}

\emph{If there exists an isometric immersion  
$f\colon M_c^n\to\Q_{\tilde c}^{n+p}$, $n\geq 2$ and $c<0$, 
with flat normal bundle  of  a complete Riemannian manifold 
$M_c^n$ of constant sectional curvature $c$ with  $c<\tilde c$, 
then $M_c^n=\Hy_c^n$}.
\vspace{1ex}

In view of Nikolayevsky's result, the following weaker 
version of the problem discussed above has already
been considered by Brander \cite{br}.

\vspace{1ex}
\noindent\emph{PROBLEM: Do isometric immersions with 
flat normal bundle of $\Hy_c^n$ into $\Q_{\tilde c}^{n+p}$ 
for $n\geq 2$ and $c<\tilde c$ exist?}
\vspace{1ex}

In this paper, we analyze the behavior of the second fundamental 
form of a possible submanifold as in the problem above, and 
conclude that it must have exponential growth, as defined next. 
\vspace{1ex}

Let $f\colon M^n\to \Q^{n+p}_{\tilde c}$ be an isometric 
immersion of a complete non-compact Riemannian manifold $M^n$. 
It is said that the second fundamental form 
$\a_f\colon TM\times TM\to N_fM$ of $f$ 
has \textit{exponential growth} if there exist $x_0\in M^n$ 
and positive constants $k,\ell\in\R$ such that 
$$
\max\left\{\|\a_f(x)\|:x\in D_r(x_0)\right\}\geq k e^{\ell r} 
$$
for any $r>r_0$ for some $r_0>0$, 
where $D_r(x_0)$ denotes the closed geodesic ball of $M^n$ of 
radius $r$ centered at $x_0$ and $\|\a_f\|$ is the norm of the
second fundamental form given by
$$
\|\a_f(x)\|^2=\sum_{i,j}\|\a_f(X_i,X_j)(x)\|^2
$$
where $X_1,\dots,X_n\in T_xM$ is an orthonormal basis.

\begin{theorem}\label{main0} If a complete $n$-dimensional Riemannian 
manifold $M^n_c$, $n\geq 2$ and $c<0$, admits an isometric immersion 
$f\colon M^n_c\to\Q_{\tilde c}^{n+p}$, $c<\tilde c$, with flat normal 
bundle then $M^n_c=\Hy^n_c$ and the second fundamental form of $f$ has 
exponential growth.
\end{theorem}

The above gives as corollary the result due to Bolotov \cite{bol} 
that there is no isometric immersion of $\Hy^n_c$ into $\R^{n+p}$ 
with mean curvature vector field of bounded length. 
\vspace{1ex}

The conclusion of Theorem \ref{main0} does not hold if the assumption
of having flat normal bundle is dropped. For instance, it was shown
by Aminov \cite{am} that the example constructed by Rozendorn
of an isometric immersion of $\Hy^2$ in $\R^5$ has no flat 
normal bundle and that the norm of its second fundamental form 
is globally bounded.
\vspace{1ex}

The aforementioned result for codimension $p=n-1$ due to Cartan has 
the following immediate consequence: 

\begin{corollary}\label{main} If there exists an isometric immersion 
$f\colon\Hy^n_c\to\Q_{\tilde c}^{2n-1}$ with $c<\tilde c$ then the 
second fundamental form of $f$ has exponential growth.
\end{corollary}

\section{The proof}

Let $f\colon M^n\rightarrow \Q^{n+p}_{\tilde c}$ 
be an isometric immersion of a Riemannian manifold $M^n$ 
into the space form $\Q^{n+p}_{\tilde c}$. 
If the immersion $f$ has flat normal bundle, 
that is, if at any point the curvature tensor of the normal 
connection vanishes, then it is a standard fact 
(cf.\ \cite{dt}) that at any point $x\in M^n$ there 
exists a set of unique pairwise distinct normal vectors 
$\eta_i(x)\in N_fM(x),\ 1\leq i\leq s(x)$, 
called the \emph{principal normals} of $f$ at $x$, and 
an associate orthogonal splitting of the tangent space as
$$
T_xM=E_{\eta_1}(x)\oplus\cdots\oplus E_{\eta_s}(x),
$$
where 
$$
E_{\eta_i}(x)=\big\{X\in T_xM:\a_f(X,Y)=\<X,Y\>\eta_i\;\;
\text{for all}\;\;Y\in T_xM\big\}.
$$
The \emph{multiplicity} of a principal normal 
$\eta_i\in N_fM(x)$  of $f$ at $x\in M^n$ is the dimension of 
the  tangent subspace $E_{\eta_i}(x)$.  
If $s(x)=k$ is constant on $M^n$, then the  maps 
$x\in M^n\mapsto\eta_i(x),  \ 1\leq i\leq k,$ 
are smooth vector fields, called the \emph{principal normal 
vector fields} of $f$. Moreover,  also the distributions 
$x\in M^n\mapsto E_{\eta_i}(x), 1\leq i\leq k$, are smooth. 
\vspace{1ex}

In the sequel, let $f\colon M_c^n\to \Q_{\tilde c}^{n+p}$, 
$c<\tilde c$, be an isometric immersion with flat normal bundle. 
Since $C=\tilde c-c>0$, it follows from the Gauss equation that
any principal normal has multiplicity one. Thus,
there exist exactly $n$ nonzero principal normal vector fields 
$\eta_1,\dots,\eta_n$ satisfying 
\be\label{gausseq}
\<\eta_i,\eta_j\>=C,\; 1\leq i\neq j\leq n.
\ee
If $X_i\in \Gamma(E_{\eta_i})$, $1\leq i\leq n$, is a unit local 
vector field  then the local orthonormal frame $X_1,\dots,X_n$ 
diagonalizes the second fundamental form of $f$, that is, 
$$
\a_f(X_i,X_j)=\delta_{ij}\eta_i,\; 1\leq i,j\leq n,
$$
where $\delta_{ij}$ is the Kronecker delta. Such a frame is called 
a \emph{principal} frame.

\begin{lemma}\label{lie} The following holds: 
\be\label{nn}
\n_{X_i} X_j= -\lambda_i X_j(1/\lambda_i) X_i,\ \ 1\leq i\neq j \leq n,
\ee
where $\lambda_i=1 \big/ \sqrt{\|\eta_i\|^2+C}$. 
\end{lemma}

\proof The Codazzi equation is equivalent to
\be\label{c1}
\nabla_{X_j}^\perp\eta_i=\<\nabla_{X_i} X_i,X_j\>(\eta_i-\eta_j)
\ee
and
\be\label{c2}
\<\nabla_{X_\ell}X_j,X_i\>(\eta_i-\eta_j)
=\<\nabla_{X_j} X_\ell,X_i\>(\eta_i-\eta_\ell)
\ee
for all $1\leq i\neq j\neq \ell\neq i\leq n$. 

The vectors $\eta_i-\eta_j$ and $\eta_i-\eta_\ell$, 
$1\leq i\neq j\neq\ell\neq i\leq n$ 
are linearly independent. Suppose otherwise that
$\eta_i-\eta_j=\mu(\eta_i-\eta_\ell)$. Taking the inner 
product with $\eta_i$ and using \eqref{gausseq} gives
$\|\eta_i\|^2=C<0$, a contradiction.

It now follows from \eqref{c2} that
$$
\n_{X_i} X_j=\Gamma_{ij}^iX_i+\Gamma_{ij}^jX_j,\;i\neq j,
$$
where $\Gamma_{ij}^k=\<\n_{X_i}X_j,X_k\>$. Since
$\Gamma_{ij}^j=\<\n_{X_i}X_j,X_j\>=0$,
then
$$
\n_{X_i}X_j=\Gamma_{ij}^iX_i=-\Gamma_{ii}^jX_i.
$$
On the other hand, taking the inner product of \eqref{c1} 
with $\eta_i$ and using \eqref{gausseq} is easily seen to 
give that $\Gamma_{ii}^j=\lambda_i X_j(1/\lambda_i)$,
as we wished.\qed

\begin{lemma}\label{lem}For each $x_0\in M^n_c$ there exists 
a diffeomorphism $F\colon U\to V$ from an open subset 
$U\subset\R^n$ endowed with coordinates $\{u_1,\ldots,u_n\}$ 
onto an open neigborhood   $V\subset M^n_c$ 
of $x_0$ such that the tangent frame
\be\label{princframe}
\sqrt{\|\eta_1\|^2+C} F_*(\d/\d u_1),\dots,
\sqrt{\|\eta_n\|^2+C} F_*(\d/\d u_n)
\ee
is orthonormal and principal. Moreover, if $M^n_c$ is complete 
and simply connected then $F\colon\R^n\to M^n_c$ is a diffeomorphism. 
\end{lemma}

\proof For the local existence, observe that Lemma \ref{lie} implies 
$$
[\lambda_i X_i,\lambda_jX_j]=0, 1\leq i\neq j\leq n.
$$

For the proof of the global part, we follow a similar argument as 
in the proof of Theorem $3$ in \cite{mo1} or Proposition $5.6$ in \cite{dt}.
Assume that $M^n_c$ is complete and simply connected. Set 
$Y_i=\lambda_iX_i$ and let $\varphi_i(x,t), \ x\in M^n_c,\ t\in\R$ 
be the one-parameter group of diffeomorphisms generated by $Y_i$. 
Since the vector fields $Y_i,\ 1 \leq i\leq n$, have bounded lengths,
it follows that $\varphi_i(x,t)$ is defined for all values of $x$ 
and $t$. Thus, for any $x\in M^n_c$, the map $t\mapsto\varphi_i(x,t)$ 
is the integral curve of $Y_i$ with $\varphi_i (x,0)=x$.
Let $x_0$ be a fixed point in $M^n_c$ and define a function 
$F=F_{x_0}\colon \R^n\to M^n_c$ by 
$$
F(t_1,t_2,\dots,t_n)
=\varphi_n(\varphi_{n-1}(\cdots\varphi_2(\varphi_1(x_0,t_1),t_2),\cdots),t_n).
$$
Since the Lie bracket $[Y_i,Y_j]$ vanishes the parameter groups 
$\varphi_i$ and $\varphi_j$ commute. This implies that
\be\label{eq32}
F_{x_0}(t+s)=\varphi_n(\varphi_{n-1}(\cdots \varphi_2
(\varphi_1(F_{x_0}(s),t_1),t_2),\cdots),t_n) = F_{F_{x_0}(s)}(t) 
\ee
where $t=(t_1,\dots,t_n)$ and $s=(s_1,\dots, s_n)$. 
Thus
$$
F_*(s)\d_i=\frac{d}{dt}|_{t=0}\ F(s_1,\dots,s_i+t,\dots,s_n) 
=\frac{d}{dt}|_{t=0} \ \varphi_i(F(s),t)
=Y_i(F(s)).
$$

We claim that $F$ is a covering map. Then this and that 
$M^n_c$ is simply connected yields that $F$ is a 
diffeomorphism, which gives the proof.

Given $x\in M^n_c$, 
let $\tilde B_{2\varepsilon}(0)$ be an open ball of radius 
$2\varepsilon$ centered at the origin such that 
$F_x|_{\tilde B_{2\varepsilon}(0)}$ is a diffeomorphism onto 
$B_{2\varepsilon}(x)=F_x(\tilde B_{2\varepsilon}(0))$. 
Set $\{\tilde{x}_\a\}_{\a\in A}=F^{-1}(x)$ and denote by
$\tilde B_{2\varepsilon}(\tilde x_\alpha)$
the open ball of radius $2\varepsilon$ centered at 
$\tilde x_\alpha$. Define a map 
$\phi_\alpha\colon B_{2\varepsilon}(x)\to
\tilde B_{2\varepsilon}(\tilde x_\alpha)$ by
$$
\phi_\alpha(y)=\tilde x_\alpha +F^{-1}_x(y).
$$
From \eqref{eq32} we obtain 
$$
F_{x_0}(\phi_\alpha(y))=F_{x_0}(\tilde x_\alpha+F^{-1}_x(y))
=F_{F_{x_0}(\tilde x_\alpha)}(F^{-1}_x(y))=F_x (F^{-1}_x(y))=y 
$$
for all $y\in B_{2\varepsilon}(x)$. Thus $F_{x_0}$ is a 
diffeomorphism from $\tilde B_{2\varepsilon}(\tilde x_\alpha)$ 
onto $B_{2\varepsilon}(x)$ having $\phi_\alpha$ as its 
inverse. In particular, this implies that 
$\tilde B_{2\varepsilon}(\tilde x_\alpha)$ and 
$\tilde B_{2\varepsilon}(\tilde x_\beta)$ are disjoint if 
$\alpha,\beta\in A$ are distinct indices. Finally, it remains 
to check that if $\tilde y \in F_{x_0}^{-1}(B_{\varepsilon}(x))$,
then $\tilde y\in \tilde B_{\varepsilon}(\tilde x_\alpha)$ 
for some $\alpha\in A$. This follows from the fact that
$$
F_{x_0}(\tilde y- F^{-1}_x(F_{x_0}(\tilde y)))
=F_{F_{x_0}(\tilde y)}(-F_x^{-1}(F_{x_0}(\tilde y)))=x.
$$
For the last equality, observe from \eqref{eq32} that for all 
$x,y\in M^n_c$ we have $F_x(t)=y$ if and only if $F_y(-t)=x$. 
\vspace{1ex}\qed

The \emph{third fundamental form ${\rm III}_f(x)$} of $f$
at $x\in M^n$ is given by 
$$
{\rm III}_f(X,Y)(x)
=\trace\<\a_f(X,\ \cdot \ ),\a_f(Y,\ \cdot \ )\>,\;\; X,Y\in T_xM.
$$
Since $\a_f$ has no kernel (that is, positive index of relative nullity), 
then ${\rm III}_f(x)$ is a positive definite inner product.

\begin{lemma}\label{flat}
The Riemannian metric ${\rm g}^0=C{\rm g}+{\rm III}_f$ is flat 
where ${\rm g}$ is the metric of $M^n_c$. Moreover, the metric 
${\rm g}^0$ is complete if {\rm g} is complete.
\end{lemma}

\proof In terms of the system of principal coordinates 
$\{u_1,\ldots,u_n\}$ given by Lemma \ref{lem}, we have 
$$
{\rm g}^0_{ij}
=C{\rm g}_{ij}+{\rm III}_f(\d/\d u_i,\d/\d u_j)
=\frac{C}{\|\eta_i\|^2+C}\delta_{ij}
+\frac{\|\eta_i\|^2}{\|\eta_i\|^2+C}\delta_{ij}
=\delta_{ij}.  
$$
Moreover, the metric ${\rm g}^0$ is complete since 
${\rm g}^0_{ij}>C{\rm g}_{ij}$.
\vspace{2ex}\qed

\noindent \emph{Proof of Theorem \ref{main0}.} By Nikolayevsky's 
result we have that $M^n_c=\Hy^n_c$. Let $F\colon \R^n\to\Hy^n_c$ be the 
global diffeomorphism given by Lemma \ref{lem}. We endow $\R^n$ with 
the pullbacks of the two metrics considered in Lemma \ref{flat} 
that are still denoted by ${\rm g}$ and ${\rm g}^0$. Notice that
$(\R^n,{\rm g}^0)$ is the standard flat Euclidean space. 

Given a smooth curve $\gamma\colon [a,b]\subset\R\to \R^n$ set 
$$
\hat S(\gamma)=\max_{t\in[a,b]}\|\a_f\|^2(F(\gamma(t))).
$$
We have from \eqref{princframe} that
$$
{\rm g}_{ij}=\frac{1}{\|\eta_i\|^2+C}\delta_{ij}
\geq\frac{1}{\hat S(\gamma)+C}\delta_{ij}.
$$ 
Then, the lengths of $\gamma$ satisfy 
\be\label{le}
L_{{\rm g}^0}(\gamma)<(\hat S(\gamma)+C)^{1/2}L_{\rm g}(\gamma).
\ee

Let $\gamma\colon [a,b]\to\R^n$ and $\tilde\gamma\colon [a,b]\to\R^n$ 
be the unique Euclidean and hyperbolic geodesics, respectively,  
joining $\gamma(a)=\tilde\gamma(a)$ to $\gamma(b)=\tilde\gamma(b)$.
From \eqref{le} we have
$$
L_{{\rm g}^0}(\gamma)\leq L_{{\rm g}^0}(\tilde \gamma)
<(\hat S({\tilde\gamma})+C)^{1/2}L_{\rm g}(\tilde \gamma).
$$
Thus, if $\gamma_{x,y}$ is the unique hyperbolic geodesic joining 
$x\neq y\in\R^n$, then the distances with respect to ${\rm g}^0$ 
and ${\rm g}$ satisfy
\be\label{di}
d_{{\rm g}^0}(x,y)<(\hat{S}({\gamma_{x,y}})+C)^{1/2}d_{\rm g}(x,y).
\ee

Fix $x_0\in\R^n$ and let $D^{\rm g}_r(x_0)$ and $D^{{\rm g}^0}_r(x_0)$ 
be the closed geodesic balls of radius $r>0$ centered at $x_0$ 
with respect to ${\rm g}$ and ${\rm g}^0$, respectively. 

It holds that
\be\label{balls}
D^{{\rm g}}_r(x_0)\subset{\rm int}
\left(D^{{\rm g}^0}_{\psi(r)}(x_0)\right), 
\ee
where 
$$
\psi(r)= r(S(r)+C)^{1/2}\;\;\mbox{and}\;\;
S(r)=\max_{x\in D^{\rm g}_r(x_0)}(\|\a_f\|^2(F(x))).
$$
In fact, if $y\in D^{{\rm g}}_r(x_0)$ we have using \eqref{di} that
$$
d_{{\rm g}^0}(x_0,y)<(\hat S(\gamma_{x_0,y})+C)^{1/2}d_{{\rm g}}(x_0,y)
\leq r (\hat S(\gamma_{x_0,y})+C)^{1/2}\leq \psi(r). 
$$
Then, we obtain using \eqref{balls} that the volumes of the geodesic 
balls satisfy 
\begin{align*}
\vol_{{\rm g}}&(D^{{\rm g}}_r(x_0))
\leq\vol_{{\rm g}}\big(D^{{\rm g}^0}_{\psi(r)}(x_0)\big)\\
&=\int_{D^{{\rm g}^0}_{\psi(r)}(x_0)} 
\Pi_{i=1}^n(\|\eta_i\|^2+C)^{-1/2} du_1\wedge\cdots\wedge du_n\\
&<\int_{D^{{\rm g}^0}_{\psi(r)}(x_0)} C^{-n/2} du_1 \wedge\cdots\wedge du_n
=C^{-n/2} \vol_{{\rm g}^0}\big(D^{{\rm g}^0}_{\psi(r)}(x_0)\big)\\  
&=r^n\left(1 + S(r)/C\right)^{n/2}\omega_n,
\end{align*}
where $\omega_n$ is the volume of the Euclidean unit $n$-ball. Since 
$\vol_{{\rm g}}(D^{{\rm g}}_r(x_0))$ is well known to grow exponentially 
with $r$ (for instance, see \cite{ch}), it follows that also $S(r)$  
grows exponentially with $r$, and thus the second fundamental form 
of $f$ has exponential growth.\qed

\begin{remark}\label{rem}{\em
It is worth mentioning that it was shown in \cite{dt95} that 
there is no isometric immersion with flat normal bundle of a complete 
Riemannian manifold $M^n_c$, $c>0$, into $\Q^{n+p}_{\tilde c}$ with 
$c<\tilde{c}$. Notice that this follows using Lemma \ref{lem}. 
}\end{remark}

\noindent Marcos Dajczer\\
IMPA -- Estrada Dona Castorina, 110\\
22460--320, Rio de Janeiro -- Brazil\\
e-mail: marcos@impa.br
\bigskip

\noindent Christos-Raent Onti\\
Department of Mathematics and Statistics\\
University of Cyprus\\
1678, Nicosia -- Cyprus
\bigskip

\noindent Theodoros Vlachos\\
University of Ioannina \\
Department of Mathematics\\
Ioannina -- Greece\\
e-mail: tvlachos@uoi.gr
 
\end{document}